\documentclass[reqno,b5paper]{amsart}
\usepackage{latexsym,amsmath,amsfonts,amssymb}
\usepackage{graphics}
\usepackage{epsfig}
\usepackage{amsmath,amscd}
\usepackage[all,cmtip]{xy}

\newtheorem{thm}{Theorem}[section]
\newtheorem{deff}[thm]{Definition}
\newtheorem{lem}[thm]{Lemma}

\newtheorem{rem}[thm]{Remark}
\newtheorem{prop}[thm]{Proposition}
\newtheorem{cor}[thm]{Corollary}

\newcommand{\vG}{\varGamma}

\newcommand{\ov}{\overline}

\def\N{{\mathbb N}}
\def\R{{\mathbb R}}

\def\eps{\varepsilon}
\def\be{\begin{equation}}

\def\ee{\end{equation}}
\def\ba*{\begin{eqnarray*}}

\def\ea*{\end{eqnarray*}}

\begin{document}

\title[Set valued integrability]{Set valued integrability in non separable Fr\'{e}chet spaces and applications}
\author[L. DI PIAZZA, V. MARRAFFA and B.
SATCO]{L. DI PIAZZA, V. MARRAFFA and B.
SATCO}
\thanks{The  authors were partially supported
 by the grant  of MURST of Italy}

\newcommand{\acr}{\newline\indent}
\address{\llap{*\,}Department of Mathematics\acr University of
Palermo\acr via Archirafi 34\acr 90123 Palermo\acr ITALY}
\email{luisa.dipiazza@unipa.it}

\address{\llap{*\,}Department of Mathematics\acr University of
Palermo\acr via Archirafi 34\acr 90123 Palermo\acr ITALY}
\email{valeria.marraffa@unipa.it}

\address{\llap{*\,}Faculty of Electrical Engineering and Computer
Science\acr "Stefan cel Mare" University of Suceava\acr Universitatii
13\acr  Suceava\acr ROMANIA}
\email{bisatco@eed.usv.ro}


\bigskip

\begin{abstract}
\noindent
We focus on measurability and integrability for set valued functions in non-necessarily separable Fr\'{e}chet spaces. We prove some properties concerning the equivalence between different classes of measurable multifunctions. We also provide useful characterizations of Pettis set-valued integrability in the announced framework. Finally, we indicate applications to Volterra integral inclusions.
\end{abstract}

\maketitle

\section{Introduction}
Due to tremendous interest in applications (e.g. control theory, optimization or mathematical economics), a wide theory was developed for the set-valued integrability in separable Banach spaces (\cite{castaing valadier}, \cite{hp} and reference inside, \cite{dp1}, \cite{dp2},
\cite{eh}, \cite{amart},  \cite{amart1}, \cite{ars1}, \cite{Z2}).
Besides, as seen in \cite{ckr1}, a similar theory for non-separable case is necessary; several authors studied the set-valued integrability without separability assumptions in
Banach spaces (see \cite{boccuto}, \cite{ckr, cr}, \cite{dp3} and \cite{musial}).\\ On the other side, it is well known that evolution partial differential equations in finite- or infinite-dimensional case can be treated as ordinary differential equations in infinite dimensional non-normable locally convex spaces. Also, the study of various differential problems on infinite intervals or that of infinite systems of differential equations naturally lead to the framework of Fr\'{e}chet spaces (as in \cite{frigon oregan}, \cite{oregan}, \cite{graef} or \cite{chang}, to cite only a few). Therefore, the study of ordinary differential equations (and inclusions) in general locally convex spaces has serious motivations.\\ But the way in this direction is paved with difficulties; indeed, even in the single-valued case, as described in \cite[Theorem 1]{lobanov sm}, the classical Peano Theorem  fails in all infinite-dimensional Fr\'{e}chet spaces. Also, a straightforward generalization of Picard-Lindel\"{o}f Theorem to Fr\'{e}chet spaces fails; existence results were obtained only under very strong assumptions on the function governing the equation (we refer to \cite{lobanov sm}, \cite{lobanov}, \cite{polewczak} or \cite{galanis3}).\\ As for the set-valued case, the situation is even more complicated since the literature concerning set-valued measure and integration in this setting is quite poor (\cite{castaing valadier, GT2}, \cite{ars}, \cite{ars2}, \cite{precup}). Differential inclusions in Fr\'{e}chet spaces have been studied, as far as we know, only in the separable case (see \cite{castaing valadier} and the references therein). One can although mention several works on differential inclusions involving at some point Fr\'{e}chet spaces, e.g. \cite{oregan} or \cite{graef}, and also the papers \cite{galanis1}, \cite{galanis2} where set differential equations using the notion of Hukuhara derivative have been considered. In all previously mentioned papers, such a space is viewed as a projective limit of a sequence of Banach spaces.\\
In the present work we focus on measurability and integrability for set valued functions in non-necessarily separable Fr\'{e}chet spaces. There are many notions of set-valued measurability and integrability available in locally convex spaces  (see \cite{barbati hess} , \cite{castaing valadier},  \cite{himmel, hp}, \cite{ars});   while, without separability, the connections between them are not always clearly stated. Therefore, in   Section 3 and in Section 4 we obtain some results concerning relations between various concepts of measurability and integrability for multifunctions defined on a probability space with non empty closed, bounded, convex
values of a Fr\'{e}chet space (in view of applications, we will consider only multifunctions defined in $[0,1]$, but all results hold true also in case of multifunctions defined in a complete probability space).

We obtain   equivalence between different classes of measurable multifunctions (see Theorem \ref{t2.4}). In these   framework we prove a version of the Pettis Measurability Theorem for set-valued functions, that, as far as we know, is new also in case the target space is a Banach space (see Theorem \ref{separable}). We also provide useful characterizations of Pettis set-valued integrability  (see Theorem \ref{p2.9} and Proposition \ref{pp}).\\
Afterwards, in Section 5  we consider applications to integral inclusions  and provide an existence result via  Kakutani fixed point Theorem for the Volterra integral inclusion
\be\label{volterra}\qquad \qquad x(t)\in x_0+ \int_0^t k(t,s)F(s,x(s))ds,  t\in [0,1].
\ee

\marginpar{}

\noindent

\marginpar{}

\section{Preliminaries}

Let $[0,1]$ be the unit interval of the real line
equipped with the usual topology and the Lebesgue measure $\mu$ and let denote by $\mathcal{L}$ the collection of all measurable sets in $[0,1]$. Throughout  this paper  $X$ is a Fr\'{e}chet space, i.e. a metrizable, complete, locally convex space. It is well known that there exists a sequence
$(p_i)_{i \in \mathbb{N}}$ of seminorms which is sufficient (that is, $p_i(x)=0$ for all $i \in \mathbb{N}$ implies that $x=0$) and increasing  such that the metric
$$
d(x,y)=\sum_{i=1}^{\infty} \frac{1}{2^i}\frac{p_i(x-y)}{1+p_i(x-y)}
$$
is invariant to translations and generates an equivalent topology on $X$.

By $X^*$ we denote the topological dual of $X$. The symbol $c(X)$ stands for the family of all nonempty closed convex subsets of $X$, while its subfamilies consisting of bounded (resp. compact) subsets will be denoted by $cb(X)$ (resp. $ck(X)$).

\noindent We consider on $c(X)$ the   Minkowski addition ($A\oplus{B}:\,=\ov{\{a+b:a\in A,\,b\in B\}}$) and the
   standard multiplication by scalars. Our unexplained terminology related to Fr\'{e}chet spaces can be found in the monographs \cite{kt} or \cite{schaeffer}.\\
Let $H$ be the corresponding Hausdorff metric on $cb(X)$, i.e.
$$
H(A,B)=\max(e_d(A,B),e_d(B,A))
$$
where the excess $e_d(A,B)$ of the set $A$ over the set $B$ is
defined as
$$
e_d(A,B)=\sup\{ d(a,B):a\in A \}=\sup   \{  \inf_{b\in B}d(a,b):
a\in A\}.
$$

In a similar way (see \cite{castaing valadier} or \cite{ars}),
for each seminorm $p_i$ we can define the semimetric
$$
H_i(A,B)=\max(e_i(A,B),e_i(B,A))
$$
where the excess $e_i(A,B)$ of the set $A$ over the set $B$ is
defined as
$$
e_i(A,B)=\sup\{ p_i(a,B),a\in A \}=\sup \{  \inf_{b\in B}p_i(a-b):
a\in A\}.$$

The families $cb(X)$ and $ck(X)$ endowed with the Hausdorff metric are complete metric spaces and the family of semimetrics  $(H_{i})_{i \in \mathbb{N}}$ generates the Hausdorff metric $H$ (see \cite{castaing touzani valadier}). If $A \subseteq X$, by the symbol $A^o$ we denote the {\it polar} of $A$, i.e. the set $\{x^* \in X^* : |\langle x^*, x \rangle| \leq 1 \ {\textrm {for \ all} }  \  x \in A\}$.

For every $C \in c(X)$ the {\it support function of} $C$ is
denoted by $\sigma( \cdot, C)$ and defined on $X^*$ by $\sigma(x^*, C) =
\sup \{ \langle x^*,x \rangle : \ x \in C\}$, for each $x^* \in
X^*$. For each $ i \in \mathbb{N}$, let $U_i$ be  the set $U_i:= \{ x \in X: p_i(x) \leq 1\}$. It is useful to recall the H${\rm \ddot{o}}$rmander equality (see \cite[Theorem II-18]{castaing valadier}):
\begin{equation}\label{Hormander}
H_i(A,B)= \sup \{ |\sigma( x^*, A)- \sigma( x^*, B)| : x^* \in U^o_i\}
\end{equation}
\noindent for all $A, B \subseteq c(X)$.

Any map $\vG: [0,1] \to c(X)$ is called a {\it multifunction}.  A function $f:[0,1] \to X$ is called a {\it selection of} $\vG$ if $$f(t) \in \vG(t)  \;  {\text{\rm a.e \ in \ }}  [0,1].$$

\noindent A multifunction $\Gamma:[0,1]\to cb(X)$ is a {\it simple multifunction} if
$\Gamma=\sum_{i=1}^p \chi_{A_i}C_i$ where $A_i \in {\mathcal{L}}$ and $C_i \in cb(X)$ are pairwise disjoint.
\noindent According to \cite{castaing touzani valadier} for every $E \in \mathcal{L}$ we define  $\int_E \Gamma(t) dt:=\sum_{i=1}^p \mu(A_i \cap E)C_i$.

\section{Measurability for multifunctions}
In this section, we focus on measurability of
multifunctions. The measurability has been defined in many different ways
in the literature concerned with the set-valued theory. We will
consider only several of these notions and try to clarify the
relationships among them.

\begin{deff}{\rm Let  $\Gamma:[0,1]\to cb(X)$ be a multifunction.
\begin{itemize}
\item[$i)$]  $\Gamma$ is said to be {\it Bochner
measurable}  (see \cite{himmel}) if it is Bochner-measurable when seen
as a single-valued function having values in the space $(cb(X),H)$
i.e. if  there exists a sequence of simple multifunctions
$\Gamma_n:[0,1]\to cb(X)$ such that
$$
\lim_{n\to\infty}H(\Gamma_n(t),\Gamma(t))=0  \ {\textrm {for \ almost \ all}} \; t\in [0,1]\;
$$
\item[$ii)$] $\Gamma$ is said to be {\it totally measurable} (see \cite{ars}) if there
exists a sequence of simple multifunctions $\Gamma_n:[0,1]\to cb(X)$
such that for every $i\in \mathbb{N}$:
\begin{equation}\label{totm}
\lim_{n\to\infty}H_i(\Gamma_n(t),\Gamma(t))=0 \ {\textrm{ for \ almost \ all}} \; t\in [0,1] \;
\end{equation}
\item[$iii)$] $\Gamma$
is said to be {\it measurable by seminorm}  or {\it p-measurable} (as in \cite{ars} or \cite{ars1}), if for
every $i\in \mathbb{N}$ there exists a sequence of simple
multifunctions $\Gamma^i_n:[0,1]\to cb(X)$ such that
$$
\lim_{n\to\infty}H_i(\Gamma^i_n(t),\Gamma(t))=0 \ {\mathrm {for \ almost \ all}} \; t\in [0,1] \;$$
\item[$iv)$] $\Gamma$ is said to be {\it Borel-measurable} (see \cite{himmel}) if it is Borel-measurable when seen as a single-valued function having values in the space $(cb(X),H)$.
\item[$v)$] $\Gamma$ is said to be {\it lower-measurable} (in \cite{himmel}) (or {\it measurable }(in \cite{hp}) or {\it Effros measurable} (in \cite{barbati hess}) if for every open subset $V$ of $X$,
$$
\Gamma^-(V)=\{t\in [0,1]: \Gamma(t)\cap V\neq \emptyset \}
$$ is measurable.
\item[$vi)$] $\Gamma$ is said to be {\it upper-measurable} (see \cite{himmel}), or  {\it strongly-measurable} (see \cite{hp})) if for every closed subset $D$ of $X$,
$$
\Gamma^-(D)=\{t\in [0,1]: \Gamma(t)\cap D\neq \emptyset \}
$$ is measurable.
\item[$vii)$] $\Gamma$ is said to be {\it scalarly-measurable} (see \cite{castaing valadier}) if
for every $x^*\in X^*$, the support functional $t\to
\sigma(x^*,\Gamma(t))$ is
measurable.
\end{itemize}}
\end{deff}

\begin{rem}\label{r0}  {\rm Observe that if $\Gamma$ is measurable by seminorm, then for each $i$ and $n \in \mathbb{N}$, the function $H_i(\Gamma_n^i(\cdot),\Gamma(\cdot)):[0,1] \rightarrow \R$ is measurable. Indeed fixed $\overline{n} \in \N$, we get that
$$ \left | H_i(\Gamma_{\overline {n}}^i(t),\Gamma(t)) - H_i(\Gamma_{\overline {n}}^i(t),\Gamma_j^i(t))\right| \leq H_i(\Gamma_j^i(t),\Gamma(t))$$
\noindent and the measurability of $H_i(\Gamma_{\overline{n}}^i(\cdot),\Gamma(\cdot))$ follows from the fact that $H_i(\Gamma_{\overline {n}}^i(\cdot),\Gamma_j^i(\cdot))$ are measurable being $\Gamma_j^i(\cdot)$ simple multifunctions for each $i$ and $j$ and $\lim_{j\to\infty}H_i(\Gamma_j^i(t),\Gamma(t))=0$.}
\end{rem}

\begin{rem}{\rm We recall that when the multifunction is a function and  $X$ is a general locally convex space, then the   Bochner measurability  implies the measurability by seminorm, while the converse implication  is true if the topology is generated by a countable family of seminorms. In particular the two concepts are the same in Fr\'echet  spaces  (see \cite[p. 247]{Garnir}).  Also, the measurability by seminorm  implies the scalar measurability  (see \cite[p. 237]{Garnir}), while according to Pettis measurability Theorem in locally convex spaces (\cite[Theorem 2.2]{blondia} and \cite[p. 248]{Garnir}),  the scalar measurability implies the measurability by seminorms if the range of the function is separable for seminorm.
Concerning the Borel measurability, simple functions are Borel measurable. However  Bochner measurable functions need not be Borel measurable in general, unless $X$ is metrizable (see \cite{Chi}), while if a function $f$ is measurable by seminorm then  the inverse images of semiballs are measurable (see \cite[p. 242]{Garnir}).

In the more general setting of multifunctions, by definition, every totally measurable multifunction is measurable by seminorm.

Moreover, it is known that
\begin{itemize}
  \item[$i)$]For multifunctions defined on   $[0,1]$ and, more generally,  on a  measure space with the disjoint hereditary additive property, the  Bochner measurability coincides with the  Borel measurability (see
      \cite{himmel} Lemma 2.5 and  p. 126).
 \item[$ii)$]  Every Borel measurable multifunction is lower measurable (see \cite[p. 131]{himmel}).
  \item[$iii)$] Every upper-measurable multifunction is lower measurable (see \cite[Proposition III.11]{castaing valadier}); the reciprocal holds for compact-valued multifunctions, as shown in \cite[Proposition III.12]{castaing valadier}.
 \item[$iv)$] The Borel measurability in general does not imply the upper-measurability (\cite[p. 131]{himmel}).
 \item[$v)$] As a consequence of
 \cite[Proposition 2.32]{hp}, every lower measurable  multifunction is also scalarly measurable.
 \item[$vi)$] The measurability by seminorm is stronger than the scalar measurability (see Remark \ref{r0} and the properties of the semimetrics $H_i$ given  in  \cite[p. 49]{castaing valadier}).
\end{itemize}

Remark also that
\begin{itemize}
 \item[$vii)$]  In separable metric spaces there are additional links between the notions presented before. In fact, Theorem III.30 in \cite{castaing valadier} shows that the lower-measurability and the upper measurability are equivalent to the fact that for every $x\in X$, the real function $t\to d(x,\Gamma(t))$ is measurable and also to the fact that the graph is measurable; Theorem III.2 in \cite{castaing valadier} states that for compact-valued multifunctions, the Bochner measurability, the lower measurability and the upper measurability coincide; besides, by Theorem III.15 in \cite{castaing valadier}, for multifunctions with convex compact values, the lower measurability coincides with the scalar measurability.
 \item[$viii)$] For Banach space-valued multifunctions without separability hypothesis, the measurability of multifunctions was investigated in \cite{cascales kadets rodriguez1}, \cite{ckr1}, \cite{cascales kadets rodriguez4}; also, the lower measurability and the scalar measurability are studied in \cite{barbati hess}; in particular, an example is provided to show that the lower measurability is strictly stronger that the scalar measurability.
 \item[$ix)$]  In \cite{hansell} some other notions of measurability of Bochner-type (via simple multifunctions) were discussed, with respect to different topologies on the hyperspace of closed sets; in particular, the Vietoris topology was considered.
\end{itemize}}
\end{rem}
\hfill$\Box$

 In order to go further into the study of various types of measurability for multifunctions in non-necessarily separable Fr\'{e}chet spaces, we need the following auxiliary result.

 \begin{lem}\label{koe} {\rm (see \cite[p. 206]{kt})}
 Let $(\alpha_i)_{i\in \mathbb{N}}$ be an increasing sequence of positive numbers and let $$\alpha=\sum_{i=1}^{\infty}\frac{1}{2^i}\frac{\alpha_i}{1+\alpha_i}.$$ Then:\\
   i) If for some $k\in \mathbb{N}$, $\alpha_k< \frac{1}{2^k}$, then
 $
 \alpha<\frac{1}{2^{k-1}}
 $.\\
 ii) If for some $m,k\in \mathbb{N}$, $\alpha<\frac{1}{2^{m+k+1}}$, then
 $
 \alpha_m<\frac{1}{2^k}
 $.
 \end{lem}

 We prove now that for multifunctions in Fr\'{e}chet spaces, the Bochner measurability, the totally  measurability and the measurability by seminorm are equivalent, as in case of Fr\'{e}chet valued functions.

\begin{thm}\label{t2.4}
Let $\Gamma:[0,1]\to cb(X)$ be a multifunction.
Then the following are equivalent:
\begin{itemize}
  \item[$i)$]  $\Gamma$ is Bochner measurable;

  \item[$ii)$]  $\Gamma$ is totally measurable;

  \item[$iii)$]  $\Gamma$ is measurable by seminorm.
  \end{itemize}
 \end{thm}
{\bf Proof.}
Let us begin by proving the equivalence of $i)$ and $ii)$.
Suppose $\Gamma$ is totally measurable. Then there exists a sequence of simple multifunctions
$\Gamma_n:[0,1]\to cb(X)$ such that for every $i\in \mathbb{N}$ there exists $N_i\subset [0,1]$ with $\mu(N_i)$=0 and
$$
\lim_{n\to\infty}H_i(\Gamma_n(t),\Gamma(t))=0, \forall t\in [0,1]\setminus N_i.$$
Fix $t\in [0,1]\setminus (\bigcup_{i=1}^{\infty}N_i)$ and let $k\in \mathbb{N}$. Since $
\lim_{n\to\infty}H_k(\Gamma_n(t),\Gamma(t))=0$, there exists $\overline{n}(t,k)\in \mathbb{N}$ such that
$$
H_k(\Gamma_n(t),\Gamma(t))<\frac{1}{2^k},\; \forall n>\overline{n}(t,k).
$$
Otherwise said, for every $n>\overline{n}(t,k)$,$$ e_k(\Gamma_n(t),\Gamma(t))<\frac{1}{2^k}\;\; {\textrm {and}}\;\; e_k(\Gamma(t),\Gamma_n(t))<\frac{1}{2^k}. $$
It follows that for each $x\in \Gamma_n(t)$, there exists $\overline{y} \in \Gamma(t)$ satisfying
$$
p_k(x-\overline{y})<\frac{1}{2^k}
$$
whence, by applying Lemma \ref{koe},
we get
$$
d(x,\overline{y})<\frac{1}{2^{k-1}}
$$
 and so,
 $$
 e_d(\Gamma_n(t),\Gamma(t))\leq\frac{1}{2^{k-1}}.
 $$
 Likewise,
 $$
 e_d(\Gamma(t),\Gamma_n(t))\leq\frac{1}{2^{k-1}},
 $$
 therefore
 $$
 H(\Gamma_n(t),\Gamma(t))\leq\frac{1}{2^{k-1}}, \; \forall n>\overline{n}(t,k).
 $$
 The Bochner measurability is thus proved.\\
Suppose now that $\Gamma$ is Bochner measurable. One can find a sequence of simple multifunctions
$\Gamma_n:[0,1]\to cb(X)$ such that there exists $N\subset [0,1]$ with $\mu(N)$=0 and
$$
\lim_{n\to\infty}H(\Gamma_n(t),\Gamma(t))=0, \forall t\in [0,1]\setminus N.$$
Fix now $i\in \mathbb{N}$ and $t\in [0,1]\setminus N$. For every $k\in \mathbb{N}$, there exists $\overline{n}(i,k,t)\in \mathbb{N}$ such that
$$
H(\Gamma_n(t),\Gamma(t))<\frac{1}{2^{i+k+1}}, \; \forall n>\overline{n}(i,k,t).
$$
It means that
$$ e_d(\Gamma_n(t),\Gamma(t))<\frac{1}{2^{i+k+1}}\;\; {\textrm {and}}\;\; e_d(\Gamma(t),\Gamma_n(t))<\frac{1}{2^{i+k+1}}. $$
So, for each $x\in \Gamma_n(t)$ one can find $\overline{y}\in \Gamma(t)$ with
$$
d(x,\overline{y})<\frac{1}{2^{i+k+1}}
$$
whence, by Lemma \ref{koe},
$$
p_i(x-\overline{y})<\frac{1}{2^k}
$$
and so,
$$
e_i(\Gamma_n(t),\Gamma(t))\leq\frac{1}{2^k}.
$$
Similarly,
$$
e_i(\Gamma(t),\Gamma_n(t))\leq\frac{1}{2^k},
$$
whence
$$
H_i(\Gamma_n(t),\Gamma(t))\leq\frac{1}{2^k}, \; \forall n>\overline{n}(i,k,t).
$$
Therefore
$$
H_i(\Gamma_n(t),\Gamma(t))\to 0 \; a.e.
$$
and the totally measurability is proved.
In order to show the equivalence of $i)$ and $ii)$ with $iii)$, it is enough to prove that $iii)$ implies $ii)$ . Suppose thus that $\Gamma$ is measurable by seminorm. Then for every $i\in \mathbb{N}$,  there exists a sequence of simple multifunctions
$\Gamma^i_n:[0,1]\to cb(X)$ such that
$$
\lim_{n\to\infty}H_i(\Gamma^i_n(t),\Gamma(t))=0\; \forall t\in [0,1] \ \ a.e.$$
We need to prove the existence of a sequence $(G_n)_n$ of simple multifunctions  on $[0,1]$ such that for each $i\in \mathbb{N}$,  $\lim_{n\to\infty}H_i(G_n(t),\Gamma(t))\to 0$ a.e.\\
By Remark \ref{r0}, for each $i \ \ {\textrm {and}} \ \ n \in \mathbb{N}$, $H_i(\Gamma^i_n(t),\Gamma(t))_n$ is measurable and also
$$ \frac{H_i(\Gamma^i_n(t),\Gamma(t))}{1+ H_i(\Gamma^i_n(t),\Gamma(t))} <1.$$
\noindent Fix  $i\in \mathbb{N}$.
\noindent The sequence  $(H_i(\Gamma^i_n(t),\Gamma(t)))_n$ converges to zero a.e., thus it converges in measure to zero. Therefore for each i there exist a simple multifunction $\Gamma^i_{\overline{n}(i)}$ and a measurable set $E_i \subset [0,1]$ with $\mu(E_i) < \frac{1}{2^i}$  such that
$H_i(\Gamma^i_{\overline{n}(i)},\Gamma(t)) < \frac{1}{2^i}$ for all $t \notin E_i$.

Then
$$ \int_0^1 \frac{H_i(\Gamma^i_{\overline{n}(i)}(t),\Gamma(t))}{1+ H_i(\Gamma^i_{\overline{n}(i)}(t),\Gamma(t))} dt=\int_{[0,1]\setminus E_i} \frac{H_i(\Gamma^i_{\overline{n}(i)}(t),\Gamma(t))}{1+ H_i(\Gamma^i_{\overline{n}(i)}(t),\Gamma(t))} dt+\int_{E_i} \frac{H_i(\Gamma^i_{\overline{n}(i)}(t),\Gamma(t))}{1+ H_i(\Gamma^i_{\overline{n}(i)}(t),\Gamma(t))} dt<  \frac{1}{2^{i-1}},$$
\noindent which implies that
$$\lim_{i\to\infty}\int_0^1 \frac{H_i(\Gamma^i_{\overline{n}(i)}(t),\Gamma(t))}{1+ H_i(\Gamma^i_{\overline{n}(i)}(t),\Gamma(t))}dt=0.$$

\noindent Then it is possible to find a subsequence  $(\Gamma^{i_k}_{\overline{n}({i_k})})_k$ of  $(\Gamma^i_{\overline{n}(i)})_i$ such that
$\lim_{{k} \to \infty}(H_{i_k}(\Gamma^{i_k}_{\overline{n}({i_k})}(t),\Gamma(t))=0$ a.e.

\noindent Defining  $G_k:=\Gamma^{i_k}_{\overline{n}({i_k})}$ we get a sequence of simple multifunctions such that $\lim_{k \to \infty}(H_k(G_k(t),\Gamma(t))=0$, for all $t\in [0,1]\setminus N$, with $\mu(N)=0$.

\noindent  Now  fix  $i\in \mathbb{N}$, $t \in  [0,1]\setminus N$ and $\eps >0$. We can find $\overline{k}=\overline{k}(i, \eps, t)$ such that if $k > \overline{k}$, we get

$$H_i(G_k(t),\Gamma(t)) \leq H_k(G_k(t),\Gamma(t)) < \eps.$$
 \noindent  Therefore the totally measurability follows  and the  proof is over.

\hfill$\Box$ \\

The SCHEMA below summarizes our discussion concerning the measurability of $cb(X)$-valued multifunctions defined on $[0,1]$ when $X$ is a general separable Fr\'{e}chet space:

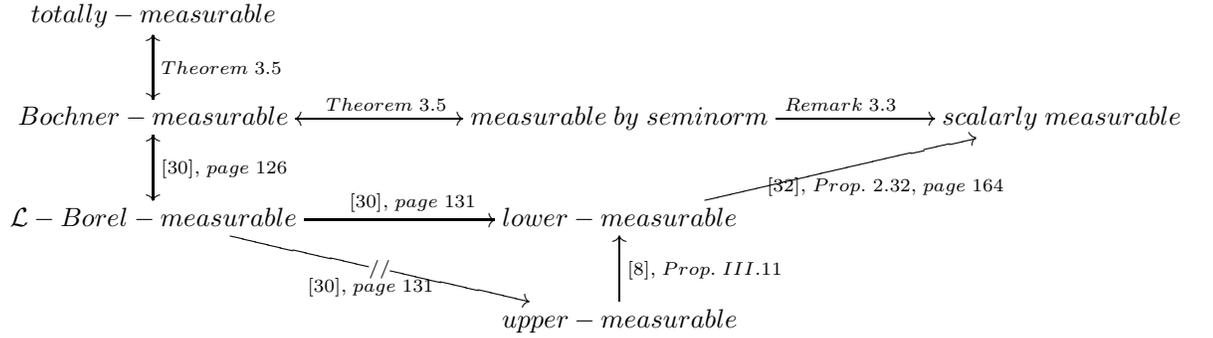
\[\xymatrixcolsep{5pc}
\xymatrix{
totally-measurable  \ar[d]^{Theorem\; 3.5}\\
Bochner-measurable \ar[u] \ar[d]^{[30],\; page \;126} \ar[r]^{Theorem\;  3.5} &measurable\; by \;seminorm \ar[l] \ar[r]^{Remark \; 3.3}  &scalarly \;measurable \\
{\mathcal L}-Borel-measurable \ar[u] \ar[r]^{\qquad [30], \;page\; 131 } \ar[rd]_{[30],\; page\; 131}|-{//} &lower-measurable \ar[ru]_{\quad [32], \;Prop.\; 2.32, \;page \;164}\\
\qquad \qquad &upper-measurable \ar[u]_{[8], \;Prop.\; III.11}}
\\
\]
By Pettis Measurability Theorem (see  \cite[Theorem 2.2]{blondia}
for its locally convex space version) it is known that any
single-valued function which is strongly measurable has almost separable
range. As far as we know, the following is the first result of
this kind in the framework of the set-valued functions taking values in Fr\'{e}chet spaces. Related to this we note only that,
in \cite[Lemma 2.5]{himmel}, it is shown that the range of a
Bochner-measurable multifunction seen like a single-valued function
having values in the space $(cb(X),H)$ is almost separable in that
space.

\begin{thm}\label{separable}
Let $\Gamma:[0,1] \rightarrow ck(X)$ be Bochner measurable. Then
there is a set $N \subseteq [0,1]$ with $\mu(N)=0$
such that $\bigcup_{t \in [0,1] \setminus N}\Gamma (t)
\subseteq X$ is separable.
\end{thm}

\noindent {\bf Proof.} Since $\Gamma$ is Bochner measurable there
is a set $N\subseteq [0,1]$ with $\mu(N)=0$  and a
sequence of simple multifunctions $\Gamma_n:[0,1]\to ck(X)$ such
that for every  $t \in [0,1] \setminus N$
$$
\lim_{n \to \infty}H(\Gamma_n(t),\Gamma(t))=0.$$ \noindent   Remark at
the first step that for each simple multifunction $G:[0,1]\to
ck(X)$, $\bigcup_{t \in [0,1]}G (t)$ is separable. Indeed, if
$G(t)=\sum_{j=1}^{l} K_j \chi_{\Omega_j}$ where, for each $j$,  $K_j$ is
compact of $X$ and $\Omega_j \in {\mathcal{L}}$, then $\bigcup_{t \in [0,1]}G (t)= \bigcup _{j=1}^{l}K_j.$
Hence it is separable (by the separability of compacts in a
metrizable space).

\noindent At the next step,
we point out that, for each $n$,  $\bigcup_{t
\in [0,1]}\Gamma_n (t)$ is separable and then there is a sequence
$(x^n_m)_{m \in \mathbb{N}}\subseteq X$ dense in $\bigcup_{t \in
[0,1]}\Gamma_n (t)$.

\noindent We claim that $$\bigcup_{t \in [0,1] \setminus
N}\Gamma (t) \subseteq \overline{D}, \; {\textrm {where}}\; D=
\{x^n_m: \ \ n,m \in \mathbb{N}\}.$$ Indeed, let $x \in \bigcup_{t
\in [0,1] \setminus N}\Gamma (t)$. Then there is $t\in [0,1]
\setminus N$ so that $x \in \Gamma (t)$, and for each $i \in
\mathbb{N}$, there is $n_i\in \mathbb{N}$ such that
$$ H(\Gamma_{n_i}(t),\Gamma(t))< \frac{1}{i},$$
\noindent which implies that there is $x_{n_i} \in \Gamma_{n_i}(t)$
such that $d(x ,x_{n_i}) <\frac{1}{i}$. Corresponding to this
$x_{n_i}$ there is $\overline{x} \in D$, such that $d(\overline{x} ,
x_{n_i}) <\frac{1}{i}$. Thus $d(\overline{x} , x) <\frac{2}{i}$
and so, the separability is proved.
\hfill$\Box$\\

\begin{rem}\label{r3}  {\rm An important aspect when dealing with multifunctions is the
existence of selections with appropriate properties (depending on
the properties of the multifunction). In this direction, when
speaking about measurability, the mostly used result is the classical Kuratowski-Ryll-Nardzewski Theorem that ensures the existence of measurable selections for measurable multifunctions in complete metrizable spaces that are separable. Without separability assumptions, in Banach spaces the existence of measurable selectors was obtained in \cite{cascales kadets rodriguez1} and \cite{ckr1}. In non-normable non-separable spaces, we cite the works \cite{hansell 83}, \cite{magerl}, \cite{hansell} or \cite{himmel}.\\
The result that will be applied in the sequel is Theorem 2.9 in \cite{himmel} which states that: {\it in a metric space any
Bochner-measurable multifunction taking convex closed bounded values  has
Bochner-measurable selections}. In \cite{himmel} the metric is supposed to be bounded but the result is available without this assumption (see  \cite[Remark 3.7]{hansell}). \\From now on $S_{\Gamma}$ will denote the family of all Bochner-measurable selections of $\Gamma$.}
\end{rem}

\section{Integrability for multifunctions}

\noindent In what follows, we are going to consider the integrability matter. In view of the applications from now on we will consider only multifunctions taking values in $ck(X)$. We observe that propositions 4.2 and 4.3 hold true also in case of multifunctions taking values in $cb(X)$.  Let us recall the following definitions:

\begin{deff}\label{int}{\rm Let  $\Gamma:[0,1]\to ck(X)$ be a multifunction.

\begin{itemize}
\item[$j)$]   $\Gamma$  is said to be {\it integrable}  (see \cite{ars}) if it is totally measurable and there exists a sequence of simple multifunctions
$\Gamma_n:[0,1]\to ck(X)$  that satisfies $\lim_{n\to\infty}H_i(\Gamma_n(t),\Gamma(t))=0 \ {\textrm{ for \ almost \ all}} \; t\in [0,1] $
 and  such that for every $i\in \mathbb{N}$:
\begin{equation}\label{e1}
\lim_{n,m \to\infty} \int_0^1 H_i(\Gamma_n(t),\Gamma_m(t))dt=0.\end{equation}

\noindent Then, for each measurable $E \subset [0,1]$, $(\int_E \Gamma
_n(t)dt)_n$ is a Cauchy sequence (\cite[p. 373]{ars}) in the Hausdorff metric. Therefore there is $x_E \in ck(X)$ such that for each $i \in  \mathbb{N}$,

\begin{equation}\label{e2}\lim_{n \to\infty} H_i\left(\int_E \Gamma_n(t)dt,x_E\right)=0, \end{equation}

\item[$jj)$]    $\Gamma$  is said to be
{\it integrable by seminorm} or {\it p-integrable} (see \cite{ars}), if it is
measurable by seminorm and for every $i\in \mathbb{N}$ there exists
a sequence of simple multifunctions $\Gamma^i_n:[0,1]\to ck(X)$ such
that
\begin{itemize}
\item[$jj_a)$]  $H_i(\Gamma^i_n(t),\Gamma(t)) \in L^1([0,1])$,
\item[$jj_b)$]   $H_i(\Gamma^i_n(t),\Gamma(t))$ converges to zero in $\mu$-measure
\item[$jj_c)$] for each
measurable $E \subset [0,1]$,
$$
\lim_{n\to\infty}\int_E H_i(\Gamma^i_n(t),\Gamma(t))dt=0,$$
\end{itemize}
\noindent and there is a unique $x_E \in
ck(X)$ such that
$$
\lim_{n\to\infty}H_i\left(\int_E\Gamma^i_n(t)dt, x_E\right)=0.$$

\item[$jjj)$] $\Gamma$  is
said to be  {\it scalarly integrable} if for every $x^* \in X^*$, the function $\sigma(x^{\ast},\Gamma(\cdot))$ is Lebesgue-integrable.

\item[$jv)$]  $\Gamma$  is
said to be  {\it Pettis-integrable } in $ck(X)$ (see \cite{castaing valadier}) if  $\Gamma$ is scalarly integrable and for every  $E \in \mathcal{L}$, there exists $x_E\in ck(X)$ such that
$$
\sigma(x^{\ast},x_E)=\int_E \sigma(x^{\ast},\Gamma(t))dt,\; {\textrm{
for}}\; {\textrm {all }}\;x^{\ast}\in X^{\ast}.
$$
\end{itemize}
In the case of Pettis integrability the set $x_E$ is denoted by $(P)\int_E \Gamma(t)dt$ and is called the {\it Pettis integral over the measurable set}  $E$, while in the case of integrability or $p$-integrability will be denoted by $\int_E \Gamma(t)dt$ and called the {\it integral of} $\Gamma$ {\it over the measurable set $E$}.}
\end{deff}

\noindent If $X$ is a Banach space the definitions of  integrability  and of   $p$-integrability  coincide. We will now
prove that this is true, more generally, even in Fr\'{e}chet spaces
(extending, in this way, the result known in the single valued case,
  \cite[Theorem 2.12]{blondia}).

\begin{prop}\label{equivalence}  A multifunction  $\vG:[0,1] \rightarrow ck(X)$ is integrable if and only if
 it is p-integrable. In such a case, the two integrals coincide.
\end{prop}

\noindent {\bf Proof.} The ``only if'' part follows by the  definition. We have to prove only the ``if'' part. By Theorem \ref{t2.4}, there is a sequence
of simple multifunctions $(\Gamma_n)_n$ such that for each $i \in
\mathbb{N}$,  $ H_i(\Gamma_{n}(t),\Gamma(t))$ converges to zero. Since
$\vG$ is $p$-integrable, for each $i \in \mathbb{N}$, there is a
simple multifunction $S_i$ such that $ H_i(S_i(t),\Gamma(t)) \in
L^1([0,1])$, and for each $E \in {\mathcal L}$,
\begin{equation}\label{pint}  \int_E H_i(S_i(t),\Gamma(t)) \ dt < \frac{1} {2^i}, \end{equation}
\noindent and moreover
 $$ H_i\left(\int_E S_i(t)dt ,x_E\right) < \frac{1}{ 2^i},$$

 \noindent where $x_E \in ck(X)$ is the $p$-integral of $\Gamma$ over $E$. We claim that this sequence $(S_n)_n$ satisfies condition (\ref{e1})  of Definition \ref{int} $j)$. Fix $k$ and $\varepsilon>0$.
 There is an $l$ such that  $l>k$ and $\frac{1}{2^l} < \varepsilon$. So   for each $n,m \geq l$ we have
 \begin{equation}\label{pint1}
\begin{split}
 H_k(S_n(t), S_m(t)) & \leq H_k(S_n(t), \Gamma(t))+  H_k(\Gamma(t), S_m(t))\\
 &\leq H_n(S_n(t), \Gamma(t))+  H_m(\Gamma(t), S_m(t))
 \end{split}\end{equation}
 \noindent which implies that $H_k(S_n(t), S_m(t)) \in L^1([0,1])$ and by (\ref{pint}) and (\ref{pint1})
 \begin{equation}\label{pint2}\int_E H_k(S_n(t), S_m(t)) < \frac{1}{2^n} + \frac{1}{2^m} < 2\varepsilon.\end{equation}
 \noindent  So condition (\ref{e1}) is satisfied. We are proving now that the integrals coincide. So again fix $k$ and $\varepsilon>0$. By  \cite[Proposition 5]{ars} for each $n, m \in {\mathbb{N}}$ and for each $E \in {\mathcal L}$, we have,
 $$H_k \left(\int_E S_n(t)dt, \int_E S_m(t)dt \right) \leq \int_E H_k(S_n(t), S_m(t))dt.$$

\noindent  Let $l\in \mathbb{N}$ be such that $l>k$ and $\frac{1}{2^l} < \varepsilon$. For each $n,m \geq l$ by (\ref{pint}) and (\ref{pint1}) we infer
\begin{equation}
\begin{split}
 H_k\left(\int_E S_n(t)dt, x_E\right) & \leq H_k\left(\int_E S_n(t)dt,\int_E S_m(t)dt\right)+  H_k \left(x_E, \int_E S_m(t)dt \right)\\
 & \leq \int_E  H_k \left(S_n(t), S_m(t)\right) dt +  H_m\left(x_E, \int_E S_m(t)dt\right) < 3 \varepsilon.
 \end{split}\end{equation}
Therefore the integrals coincide.
 \hfill$\Box$

As for the relationship with the Pettis integrability, it is not difficult to see that if $\Gamma$ is integrable, then it is also Pettis
integrable.

\begin{prop}\label{intp}
Let $\Gamma:[0,1]\to ck(X)$ be integrable. Then it is also Pettis integrable in $ck(X)$ and the two integrals coincide.
\end{prop}
\noindent {\bf Proof.} By \cite[Proposition 1]{ars}, for every seminorm $p_i$, the scalar function $H_i(\Gamma(\cdot),\{0\})$ is integrable. By  (\ref{Hormander})  the scalar integrability of $\Gamma$ follows. In order to obtain the Pettis integrability, it is enough to prove that the integral of $\Gamma$ in the sense of Definition \ref{int} $j$) satisfies the equality
$$
\sigma\left(x^{\ast},\int_E \Gamma(t)dt \right)=\int_E \sigma(x^{\ast},\Gamma(t))dt,
$$
for every $E \in \mathcal{L}$ and for every $x^* \in X^*$ .  Since $\Gamma$ is Bochner measurable, there exists a sequence of simple multifunctions $\Gamma_n=\sum_{i=1}^{p_n} \chi_{A^n_i}C^n_i$ where $A^n_i \in {\mathcal{L}}$ and $C^n_i \in ck(X)$ like in Definition \ref{int} $j$). Fix  $x^* \in X^*$ and let $i$ such that $x^* \in U_i^0$, where  $U_i^0$ denotes the polar of the set $U_i$. Then taking in account that for each $i, \ n $, and for each $E \in {\mathcal{L}}$, $H_i(\int_E\Gamma(t), \int_E \Gamma_n(t)) \leq \int_E H_i(\Gamma(t), \Gamma_n(t))$ (see \cite[Proposition 5]{ars}),  by (\ref{Hormander}), (\ref{e1}) and (\ref{e2})  we get
\begin{eqnarray*}
\sigma\left(x^{\ast},\int_E \Gamma(t)dt \right)&= &\lim_{n\to\infty} \sigma\left(x^{\ast},\int_E \Gamma_n(t)dt \right)\\
&= &\lim_{n\to\infty} \sigma\left(x^{\ast}, \sum_{i=1}^{p_n} \mu(A^n_i \cap E)C^n_i \right)\\
&= &\lim_{n\to\infty} \sum_{i=1}^{p_n} \mu(A^n_i \cap E) \sigma\left(x^{\ast}, C^n_i \right)\\
&= &\lim_{n\to\infty} \int_E \sigma\left(x^{\ast},\Gamma_n(t)\right)dt\\
& =&\int_E \sigma\left(x^{\ast},\Gamma(t)\right)dt.
\end{eqnarray*}
\hfill$\Box$

\noindent In applications to differential inclusions we will need other results, which are well known for Banach valued multifunctions (see \cite{eh}) and also for single-valued Pettis integrable functions in Fr\'{e}chet spaces (\cite{ali chakraborty 1}).\\ We recall that a map $M:\mathcal{L} \rightarrow ck(X)$ is called a  {\it weak multimeasure} or simply a {\it multimeasure} if $\sigma(x^*,M(\cdot))$ is a measure, for every $x^* \in X^*$.

\noindent  In order to prove next proposition we will use the following characterization of the $ck(X)$ sets (see \cite[Section 7]{moreau}   for the general case and \cite[Proposition 1.5]{eh}  for the Banach valued case): {\it a nonempty set    $C \in c(X)$   belongs to $ck(X)$ if and only if,  for every $0$-neighborhood $U$ in $X$,  the restriction of $\sigma( \cdot, C)$  to the polar $U^o$ of $U$ is weak$^*$-continuous}.

\begin{thm}\label{p2.9} Let $\Gamma:[0,1] \rightarrow ck(X)$ be Bochner-measurable. Then $\Gamma$ is Pettis integrable in $ck(X)$ if and only if for each equicontinuous set $K \subset X^*$, the set $\{\sigma(x^*,\Gamma(\cdot)) : x^* \in K \}$ is uniformly integrable.
\end{thm}
\noindent {\bf Proof.} Assume first that the multifunction $\Gamma$ is Pettis integrable in $ck(X)$. In order to prove that for each equicontinuous set $K \subset X^*$, the set $\{\sigma(x^*,\Gamma(t)) : x^* \in K \}$ is uniformly integrable, it is enough to show that for each $i$  the set $\{\sigma(x^*,\Gamma(t)) : x^* \in U_i^o \}$ is uniformly integrable (see \cite[p. 258]{kt}). So fix $i$.
For each  $A \in {\mathcal L}$, let
$$ M(A)= (P)\int_A \Gamma(t) dt. $$
\noindent Then $M$ is a weak multimeasure with convex compact  values and  by  \cite[Proposition 3]{GT} $M$ is a normal multimeasure. Therefore   $M$ is $\sigma$-additive with respect to the uniformity of the Hausdorff distance, that is if $A= \cup_{i=1}^{\infty} A_i$, with $(A_i)_{i\in \mathbb{N}}$ pairwise disjoint, then
$$ \lim_{n \rightarrow \infty}H \left(M(A),{\sum_{j=1}^n M(A_j)} \right) =0,$$
where ${M(A) + M(B)}= M(A \cup B)$, if $A \cap B =\emptyset$.
\noindent So  by Lemma 3.3 also for each $i$
$$ \lim_{n \rightarrow \infty}H_i\left(M(A), {\sum_{j=1}^n M(A_j)}\right) =0.$$
By the H${\rm \ddot{o}}$rmander equality  (\ref{Hormander}) we have
$$ H_i\left(M(A), {\sum_{j=1}^n M(A_j)}\right)=\sup\left\{ \left|\sigma(x^*,M(A) )- \sigma\left(x^*, {\sum_{j=1}^{n} M(A_j)}\right) \right| : x^* \in U_i^o\right\}.$$
\noindent Therefore
$$ \lim_{n \rightarrow \infty}\left[\sup\left\{ \left|\sigma(x^*,M(A)) - \sigma\left(x^*,  {\sum_{j=1}^n M(A_j)}\right) \right| : x^* \in U_i^o\right\}\right]=0.$$
\noindent Since
 \begin{equation}
\begin{split}
&\lim_{n \rightarrow \infty}\left[\sup\left\{\left|\sigma(x^*,M(A)) - \sigma\left(x^*,  {\sum_{j=1}^n M(A_j)}\right)\right| : x^* \in U_i^o\right\}\right]\\
= &\lim_{n \rightarrow \infty}\left[\sup\left\{\left|\sigma\left(x^*,M\left(\bigcup_{j=n+1}^{\infty} A_j\right)\right)\right| : x^* \in U_i^o\right\}\right]\\
= & \lim_{n \rightarrow \infty}\left[\sup\left\{\left|\sigma\left(x^*,\int _{\bigcup_{j=n+1}^{\infty} A_j} \Gamma(t) dt \right)\right| : x^* \in U_i^o\right\}\right]
\end{split}
\end{equation}
\noindent it follows that

$$\lim_{n \rightarrow \infty}\left[\sup\left\{\left|\sigma\left(x^*,\int _{\bigcup_{j=n+1}^{\infty} A_j}\Gamma(t) dt \right)\right| : x^* \in U_i^o\right\}\right]=0.$$

\noindent This implies that the collection of scalar measures $\{\sigma(x^*,M(\cdot)) : x^* \in U_i^o \} $ is uniformly $\sigma$-additive. Moreover, for every $x^*\in X^* $,   $\sigma(x^*,M( \cdot))$ is a real signed measure. Now let  $M(\mathcal{L}):=\bigcup_{A\in \mathcal{L}}M(A)$  and let $B^+$ and $B^-$  $ \in {\mathcal L}$ the corresponding Hahn decomposition of $[0,1]$. Then we have
\begin{equation*}\begin{split}
\sigma(x^*,M(\mathcal{L})) &= \sigma(x^*,\cup_{A \in {\mathcal  {L}}}M(A) )= \sup_{A \in {\mathcal {L}}} \sigma (x^* , M(A))\\
&=\sigma (x^* , M(B^+)) = < x^*, {\overline x} >,
\end{split}
\end{equation*}
\noindent where $ {\overline x} \in M(B^+)$ and the last equality follows from the fact that $M(B^+) \in ck(X)$. Since $x^* \in X^* $ is arbitrary, by the James' Theorem we get that ${\overline {M({\mathcal {L}})}}$ is weakly compact. Then  $M({\mathcal{L}})$ is bounded.
 Therefore there is a constant $\beta_i$ such that, for each $x \in M(\mathcal{L})$,  $p_i(x) \leq \beta_i$. This yields that for each $x^* \in U_i^o$
$$|\sigma(x^*,M(\mathcal{L}))| = \sup \{ |x^*(x) |: x \in M(\mathcal{L})\}  \leq  \beta_i$$ and the family $\{\sigma(x^*,M(\cdot)) : x^* \in U_i^o \} $ is bounded. Moreover, for every $x^* \in U_i^o$, the measure $\sigma(x^*,M(\cdot))$ is absolutely continuous with respect to $\mu$.  It follows from \cite[Corollary I.2.5]{diestel uhl}  that the measures of the family  are also uniformly absolutely continuous with respect to $\mu$. Since  $\sigma(x^{\ast},M(A))=\int_A \sigma(x^{\ast},\Gamma(t))dt$ for every $A \in \mathcal{L}$,  the uniform integrability of the set $\{\sigma(x^*,\Gamma(t)) : x^* \in U_i^o \}$ follows.

\noindent To prove the converse, we observe that  in order to obtain  the Pettis integrability of $\Gamma$  in $ck(X)$, we have to show that there exists $C \in ck(X)$
    such that $\int_0^1 \sigma(x^*,\Gamma(t))\,dt= \sigma(x^*,C)$
  for all functionals $x^*$. We shall
    prove first that the sublinear function $a:X^*\to(-\infty,+\infty)$ given by
    $a(x^*):=\int_0^1 \sigma(x^*,\Gamma(t))\,dt$
    is $w^*$-lower semi-continuous, i.e.
    that for each real $\alpha$ the set
    $Q(\alpha):=\{x^*\in X^*: a(x^*)\leq\alpha\}$ is $w^*$-closed. By \cite[Theorem 6.4]{schaeffer}
    it suffices to show that,  for each $i$,  $Q(\alpha)\cap U^o_i$ is  $w^*$-closed.
 Since $\Gamma $ is Bochner measurable, by Theorem \ref{separable} there is a set $N \subseteq [0,1]$ with $\mu(N)=0$
such that $\bigcup_{t \in [0,1] \setminus N}\Gamma (t) \subseteq X$ is separable. Thus restricting to the closed linear span of the set $\bigcup_{t \in [0,1] \setminus N}\Gamma (t)$, by \cite[p. 259]{kt}  it is enough to prove that  $Q(\alpha)\cap U^o_i$ is sequentially $w^*$-closed.
So let $x_n^*\in Q(\alpha)\cap U^o_i$ be  such that $x_n^*\to x_0^*$ in $\sigma(X^*,X)$. Since $\Gamma(t) \in ck(X)$, applying the $w^*$-continuity
of all $\sigma(\cdot,\Gamma(t))$, we get the pointwise convergence of $ \sigma(x_n^*,\Gamma(t))$ to $\sigma(x_0^*,\Gamma(t))$. As  by the Alaoglu-Bourbaki Theorem $U^o_i $ is compact in the $w^*$-topology, then  it is also equicontinuous. So  by hypothesis  the set $\{\sigma(x^*,\Gamma(t)) : x^* \in U^o_i \}$ is uniformly integrable. Therefore by the Vitali convergence Theorem   we have
   \begin{equation}\begin{split}
   a(x_0^*) & =\int_0^1 \sigma(x_0^*,\Gamma(t))\,dt= \int_0^1 \lim_n \sigma(x_n^*,\Gamma(t))\,dt\\
    & = \lim_n \ \int_0^1 \sigma(x_n^*,\Gamma(t))\,dt=\lim_n a(x_n^*)\leq\alpha
    \end{split}
    \end{equation}

\noindent So, $a$ is $w^*$-lower semi-continuous, and according to \cite[Theorem 5]{hormander}, there exists a closed convex set $C \subset X$ such that $a(x^*)= \sigma(x^*,C)$. Moreover if $x_n^*\in  U^o_i$ is such that $x_n^*\to x_0^*$ in $\sigma(X^*,X)$, then proceeding as in the previous part, we get  that  $ a(x_0^*)  = \lim_n a(x_n^*)$. Since each polar of a $0$-neighborhood is metrizable in the  $\sigma(X^*,X)$ topology, this means that the restriction of $\sigma( \cdot, C)$  to  $U^o_i$ is weak$^*$-continuous. Then the compactness of $C$ follows and $\Gamma$ is Pettis integrable in $ck(X)$.  \hfill$\Box$

\begin{rem} {\rm We  note that, in the proof of the previous characterization of Pettis integrability, the Bochner measurability is not needed in order to prove the uniform integrability of  the set $\{\sigma(x^*,\Gamma(\cdot)) : x^* \in K \}$ when the multifunction $\Gamma(\cdot)$ is Pettis integrable.}
\end{rem}
Without the Bochner measurability hypotheses, we can prove the following characterization:
\begin{prop}\label{pp}
Let $\Gamma:[0,1] \rightarrow ck(X)$ be  scalarly  measurable. Then the following are equivalent:
\begin{itemize}
  \item[$i)$]  $\Gamma$ is Pettis integrable in $ck(X)$;

  \item[$ii)$]  for each $A \in \mathcal{L}$ and for each $0$-neighborhood $U_i$ the map

  $$\phi_A^{\Gamma}: X^*\to(-\infty,+\infty), \ \ \ \ x^* \rightarrow \int_A \sigma(x^*,\Gamma(\cdot))dt$$
  \noindent restricted to  the polar $U^o_i$ of $U_i$ is weak$^*$-continuous.

 \end{itemize}

\end{prop}
\noindent {\bf Proof.} Assume first that $\Gamma$ is Pettis integrable in $ck(X)$. Then, for each $A \in \mathcal{L}$, there exists $C_A \in ck(X)$ such that for each $x^* \in X^*$
$$ \sigma(x^*,C_A) = \int_A \sigma(x^*,\Gamma(\cdot))dt.$$
\noindent Since $C_A \in ck(X)$, for every $0$-neighborhood $U_i$ in $X$ the restriction of $\sigma( \cdot, C_A)$  to the polar $U^o_i$ of $U_i$ is weak$^*$-continuous. By previous equality the thesis follows.

\noindent Conversely assume that $ii)$ holds. We have to show that for each $A \in \mathcal{L}$, there exists $C_A \in ck(X)$ such that
$ \sigma(x^*,C_A) = \int_A \sigma(x^*,\Gamma(\cdot))dt$ for each $x^* \in X^*$.
 Fix $A \in \mathcal{L}$ and $U_i$ in $X$.
At first we want to prove that the sublinear  map $\phi_A^{\Gamma}:X^*\to(-\infty,+\infty)$
    is $w^*$-lower semi-continuous. To do this we have to show
    that, for each real $\alpha$, the set
    $Q(\alpha):=\{x^*\in X^*: \phi_A^{\Gamma}(x^*)\leq\alpha\}$ is $w^*$-closed. By \cite[Theorem 6.4]{schaeffer}
    it suffices to show that, for each $i$,  $Q(\alpha)\cap U^o_i$ is  $w^*$-closed and this  follows  by the hypothesis that the map $\phi_A^{\Gamma}$ restricted to  $U^o_i$ is weak$^*$-continuous.
Therefore, $\phi_A^{\Gamma}$ is $w^*$-lower semi-continuous, and according to \cite[Theorem 5]{hormander}, there exists a closed convex set $C_A \subset X$ such that $\phi_A^{\Gamma}(x^*)= \sigma(x^*,C_A)$. Since the restriction of $\sigma( \cdot, C_A)$  to  $U^o_i$ is weak$^*$-continuous,  the compactness of $C_A$ follows and $\Gamma$ is Pettis integrable in $ck(X)$.  \hfill$\Box$

\begin{cor}\label{cc}
Let $\Gamma:[0,1] \rightarrow ck(X)$ be Pettis integrable in $ck(X)$ and $G:[0,1] \rightarrow ck(X)$ be scalarly  measurable and satisfying the condition
$$
G(t)\subset \Gamma(t),\; \forall t \in [0,1].
$$
Then the multifunction $G$ is Pettis integrable in $ck(X)$.
\end{cor}
\noindent {\bf Proof.}  Since  we have

$$ -\sigma(-x^*, \Gamma(t)) \leq \sigma(x^*, G(t)) \leq \sigma(x^*, \Gamma(t)), $$
\noindent the integrability of $\sigma(x^*, G)$ follows. Now fix $A \in \mathcal{L}$. The map $\varphi_A^G$ is subadditive and satisfies  $\varphi_A^G(x^*) \leq \varphi_A^{\Gamma}(x^*)$, for all $x^* \in X^*$. Hence
$$-\varphi_A^{\Gamma}(y^*-x^*)\leq -\varphi_A^G(y^*-x^*)\leq \varphi_A^G(x^*)-\varphi_A^G(y^*)\leq \varphi_A^G(x^*-y^*) \leq \varphi_A^{\Gamma}(x^*-y^*).$$

Since $\Gamma$ is Pettis integrable, by Proposition \ref{pp} the map $\varphi_A^{\Gamma}$ is weak$^*$-continuous so, it is weak$^*$-continuous at the origin and then the previous inequality implies that also the map $\varphi_A^G$ is weak$^*$-continuous. Applying again Proposition \ref{pp} we get the Pettis integrability in $ck(X)$ of $G$.
\hfill$\Box$

\vspace{2ex}

By previous result we get at once

\begin{cor}\label{cc1}
Let $\Gamma:[0,1] \rightarrow ck(X)$ be Pettis integrable in $ck(X)$. Then every scalarly measurable selection of $\Gamma$ is Pettis  integrable.
\end{cor}

\vspace{2ex}

A  property similar to Corollary \ref{cc} can be proved for integrability in a stronger sense.
\begin{prop}\label{l101} Let $\Gamma:[0,1] \rightarrow ck(X)$ be integrable and $G:[0,1] \rightarrow ck(X)$ be Bochner measurable and satisfying the condition
$$
G(t)\subset \Gamma(t),\; \forall t\in [0,1].
$$
Then the multifunction $G$ is integrable.
\end{prop}
\noindent {\bf Proof.} By Theorem \ref{t2.4} it follows that $\Gamma$ is measurable by seminorm.
Since for each $i \in \mathbb{N}$
$$0 \leq H_i(G(\cdot), \{0\}) \leq H_i(\Gamma(\cdot), \{0\}),$$
\noindent and since from \cite[Proposition 1]{ars} the function $H_i(\Gamma(\cdot), \{0\})$ is integrable, by \cite[Theorem 3]{ars}, the
integrability of the multifunction $G$ follows.
\hfill$\Box$\\

When studying multivalued differential problems, a special attention must be paid to integrable selections of the multifunction on the right hand side. The results proved below are therefore  important in this framework.

\begin{prop}\label{closed}
Let $\Gamma:[0,1] \rightarrow ck(X)$ be Bochner-measurable and Pettis integrable in $ck(X)$. Then  for each $E \in \mathcal{L}$ the set
$$I_E (\Gamma):=
\left\{ (P)\int_E f(t)dt: f\; {\rm {Pettis}}\;{\rm{-integrable}}\;{\rm {selection}}\; {\rm {of}} \; \Gamma \right\}$$
\noindent is closed.
\end{prop}
\noindent {\bf Proof.}
By Corollary \ref{cc1} for any measurable $E$, the set $I_E (\Gamma)$ is nonempty. Since $\Gamma$ is Bochner-measurable, by
Theorem  \ref{separable}  there is a set $N \subseteq [0,1]$ with $\mu(N)=0$
such that $\bigcup_{t \in [0,1] \setminus N}\Gamma (t) \subseteq X$ is separable, and being  $\Gamma$   Pettis integrable, we can assume that  $\bigcup_{t \in [0,1]}\Gamma (t) \subseteq X$ is separable.
Thus restricting to the closed linear span of the set $\bigcup_{t \in [0,1] }\Gamma (t)$, by \cite[p. 259]{kt}  the proof of the closeness of the set $I_E (\Gamma)$ follows  as in \cite[Proposition 5.2]{eh} with suitable changes.
\hfill$\Box$\\

\begin{thm}\label{Pettis}
Let $\Gamma:[0,1] \rightarrow ck(X)$ be Bochner-measurable and Pettis integrable in $ck(X)$. Then  for each $E \in \mathcal{L}$
$$\int_E \Gamma(t) dt=\left\{ (P)\int_E f(t)dt: f\; {\rm {Pettis}}\;{\rm{-integrable}}\;{\rm {selection}{\; {\rm {of}} \; \Gamma .}}\right\}.$$
\end{thm}
\noindent {\bf Proof.}
By Lemma \ref{closed} for any measurable $E$, the set
$$I_E (\Gamma):=
\left\{ (P)\int_E f(t)dt: f\; {\textrm {Pettis}}\;{\textrm{-integrable}}\;{\textrm {selection}{\; {\textrm {of}} \; \Gamma .}}\right\}$$
\noindent is closed. Now we want to prove that for each $x^* \in X^*$ and for each $E \in {\mathcal{L}}$,

$$\sigma(x^*, I_E (\Gamma)) = \int_E \sigma(x^*, \Gamma(t)) dt.$$
\noindent
Let  $x^* \in X^*$ and  $f \in S_{\Gamma}$. Then  $\langle
x^*,f(t)\rangle \leq \sigma(x^*,\Gamma(t))$ for every
$t \in [0,1]$, therefore
\begin{equation*}
\int_{E}\langle
x^*,f(t) \rangle dt \leq  \int_{E} \sigma(x^*,\Gamma(t)) dt.
\end{equation*}
This gives that \begin{equation}\label{e3}
\sigma(x^*, I_E (\Gamma)) \leq \int_E \sigma(x^*, \Gamma(t)) dt.
\end{equation}
To prove the reverse inequality, let us fix $x^*$ and consider
the multifunction $G:[0,1] \to ck(X)$ defined by

$$
G(t):\,=\{x\in \Gamma(t):\langle x^*,x\rangle =
\sigma(x^*,\Gamma(t))\}.
$$
Since $\Gamma$ has  compact convex values, for each $t \in
[0,1]$, $G(t)$ is nonempty. By  \cite[Lemma 3]{valadier}, we infer
that $G$ is scalarly measurable. For each $t \in [0,1]$,
$G(t)\subset \vG(t)$,  then by Corollary \ref{cc} $G$ is Pettis integrable.
Let $g$ be a Pettis integrable selection of $G$ and then  also  of $\Gamma$.
 Clearly, for all $t\in
[0,1]$ it satisfies the equality
 $$
 \langle x^*,g(t) \rangle =
\sigma(x^*,\Gamma(t)).
$$
Therefore
$$
\int_{E} \langle x^*,g(t)\rangle dt = \int_{E}
\sigma(x^*,\Gamma(t))dt,
$$
whence
\begin{equation}\label{e4}
\sigma(x^*, I_E (\Gamma))  \geq \int_{E} \langle x^*,g(t)\rangle dt =
\int_{E}  \sigma(x^*,\Gamma(t))dt .
 \end{equation}
 By (\ref{e3}) and (\ref{e4}) the
assertion follows.
\hfill$\Box$

\section{Volterra integral inclusions in Fr\'{e}chet spaces}

In this section we apply previously obtained  results to provide existence results for the integral problem (\ref{volterra}). \\

The first one deals with the Pettis integrability notion, therefore it seems useful to recall some basic properties of the primitive of a Pettis integrable function: $t\in [0,1]\mapsto {\textrm {(P)}}\int_0^t f(s)ds$. It is, by definition, weakly continuous. But thanks to the characterization given in \cite{ali chakraborty 1} (which was generalized by us to the set-valued case in Theorem \ref{p2.9}), it is in fact continuous. Moreover, it is pseudo-differentiable in the sense described below:
\begin{deff}  {\rm A function $F:[0,1]\rightarrow X$ is said to be {\it pseudo-differentiable}
with a pseudo-derivative $f$ if, for every $x^{\ast}\in X^{\ast}$,
there exists $N(x^{\ast})\subset [0,1]$ of null measure such that
$\left\langle x^{\ast},F \right\rangle$ is differentiable on
$[0,1]\setminus N(x^{\ast})$ and its derivative is $\left\langle
x^{\ast},F(t) \right\rangle'= \left\langle x^{\ast},f(t)
\right\rangle$, for every $t \in [0,1]\setminus N(x^{\ast})$.}
\end{deff}
Note that the weak differentiability of the Pettis primitive (namely, the existence of a null-measure set $N$ independent of $x^*$ except which the property in the preceding definition holds) is not valid even in separable Banach spaces (as proved e.g. in \cite{girardi}).\\

We proceed to give the first existence result. Our proof uses a technique  similar to that applied for a particular case in \cite[Theorem VI-7]{castaing valadier}, where $X$ is separable. It is based on Kakutani-Ky Fan Theorem (\cite[Theorem 6.5.19]{papa}). We do not impose the separability to the Fr\'{e}chet space $X$, but applying Theorem  \ref{separable}, we get by this lack by considering Bochner-measurable multifunctions in a general Fr\'{e}chet space.
\begin{thm}\label{Volterra}
Let $\Gamma:[0,1]\to ck(X)$ be Bochner-measurable and Pettis-integrable in $ck(X)$ and $F:[0,1]\times X\to ck(X)$ satisfy the following assumptions:\\
1) for every $x\in C([0,1],X)$, the multifunction $t\mapsto F(t,x(t))$ is Bochner-measurable;\\
2) for every $t\in [0,1]$, $x\mapsto F(t,x)$ is upper semi-continuous;\\
3) for every $t\in [0,1]$ and $x\in X$, $F(t,x)\subset \Gamma(t)$.\\
Let $k:[0,1]\times [0,1]\to \mathbb{R}$ be such that for each $t\in [0,1]$, the function $k(t,\cdot)\in L^{\infty}([0,1])$ and $t\mapsto k(t,\cdot)$ is $L^{\infty}$-continuous. Then the integral problem (1)  has solutions $x\in C([0,1],X)$.
\end{thm}
{\bf Proof.}
Consider $$\mathcal{X}=\left\{ z:[0,1]\to X :  \;  z(t)=x_0+{\textrm {(P)}}\int_0^t k(t,s)f(s)ds, \; \forall t\in [0,1], \;  f\in S_{\Gamma} \right\}$$
where we recall that $S_{\Gamma}$ is the family of all Bochner-measurable selections and, (by Corollary \ref{cc1}, of all Pettis-integrable selections) of $\Gamma$. By Theorem 2.9 in \cite{himmel}, $S_{\Gamma}$  is non-empty and, as  $k(t,\cdot)\in L^{\infty}([0,1])$, it follows that for each $f\in S_{\Gamma}$, $k(t,\cdot)f(\cdot)$ is Pettis-integrable (\cite[Theorem 2.3]{ali chakraborty 1}) and so, $\mathcal{X}$ is non-empty as well.\\
We obtain by the particular single-valued case of Theorem \ref{p2.9} (see \cite{ali chakraborty 1}) that $\mathcal{X}\subset C([0,1],X)$ since for every seminorm $p_i$,
\begin{eqnarray*}
&&\sup_{x^*\in U_i^0} \langle x^*,z(t)-z(t_0) \rangle \\
 && \leq \sup_{x^*\in U_i^0} \int_0^t \left| \langle x^*,(k(t,s)-k(t_0,s))f(s) \rangle \right| ds +\sup_{x^*\in U_i^0} \int_{t_0}^t \left| \langle x^*,k(t_0,s)f(s)\rangle \right| ds \\
&& \leq \|k(t,\cdot)-k(t_0,\cdot)\|_{\infty} \sup_{x^*\in U_i^0} \int_{0}^{t} \left| \langle x^*,f(s) \rangle \right| ds+\sup_{t\in[0,1]}\|k(t,\cdot)\|_{\infty} \sup_{x^*\in U_i^0} \int_{t_0}^{t} \left| \langle x^*,f(s) \rangle \right| ds.
\end{eqnarray*}
We intend to prove now that $\mathcal{X}$ is compact in the topology of uniform convergence. By Ascoli's Theorem (\cite[p. 233]{kelley}), it is enough to prove that:\\
i) $\mathcal{X}$ is equi-continuous;\\
ii) for every $t\in [0,1]$, the subset $\{z(t): \ z \in \mathcal{X} \}\subset X$ is relatively compact;\\
iii) $\mathcal{X}$ is closed.\\
In order to prove the condition i), see that by Theorem \ref{p2.9} for every seminorm $p_i$ and every
$\varepsilon>0$ there exists
$\delta_{\varepsilon,p_i}>0$ such that for all $t_1 < t_2\in [0,1]$ with
$t_2-t_1<\delta_{\varepsilon,p_i}$,
$\|k(t_1,\cdot)-k(t_2,\cdot)\|_{\infty}\leq \varepsilon$ and
$$
\int_{t_1}^{t_2} |\sigma(x^{\ast},\Gamma(s))|ds<\varepsilon,\quad \forall \; x^* \in U_i^0.
$$
It follows that for all $f\in S_{\Gamma}$,
$$
\int_{t_1}^{t_2} \left|\langle x^{\ast},f(s) \rangle  \right|ds \leq \max \left( \int_{t_1}^{t_2} |\sigma(x^{\ast},\Gamma(s))|ds, \int_{t_1}^{t_2} |\sigma(-x^{\ast},\Gamma(s))|ds \right) < \varepsilon
$$
whenever
$|t_1-t_2|<\delta_{\varepsilon,p_i}$ for all $x^* \in U_i^0$. So, for every $z\in \mathcal{X}$ and every such $t_1<t_2\in [0,1]$,
\begin{eqnarray*}
&&\sup_{x^*\in U_i^0} \langle x^*,z(t_2)-z(t_1)\rangle \\
&&\leq \sup_{x^*\in U_i^0}  \int_0^{t_1} \left|\langle x^*,(k(t_2,s)-k(t_1,s))f(s)\rangle \right|ds +\sup_{x^*\in U_i^0} \int_{t_1}^{t_2} \left|\langle x^*,k(t_2,s)f(s)  \rangle \right| ds \\
&& \leq \|k(t_1,\cdot)-k(t_2,\cdot)\|_{\infty} \sup_{x^*\in U_i^0} \int_{0}^{1} \left| \langle x^*,f(s) \rangle \right| ds+\sup_{t\in[0,1]}\|k(t,\cdot)\|_{\infty} \sup_{x^*\in U_i^0} \int_{t_1}^{t_2} \left| \langle x^*,f(s) \rangle \right| ds \\
&&  \leq \varepsilon \sup_{x^*\in U_i^0} \max \left( \int_{0}^{1} |\sigma(x^{\ast},\Gamma(s))|ds, \int_{0}^{1} |\sigma(-x^{\ast},\Gamma(s))|ds \right)+ \varepsilon \sup_{t\in[0,1]}\|k(t,\cdot)\|_{\infty}
\end{eqnarray*}
and so, the equicontinuity is satisfied.\\
The condition ii) immediately comes from the fact that by Theorem \ref{Pettis}
$
\{z(t) : \ z \in \mathcal{X} \}\subset x_0+{\textrm {(P)}}\int_0^t k(t,s)\Gamma(s)ds
$
and by Theorem \ref{p2.9} which yields the Pettis-integrability of $k(t,\cdot)\Gamma(\cdot)$.\\
As for the condition iii), take a sequence $(z_{n})_{n}\subset \mathcal{X}$ convergent to $z\in C([0,1],X)$ and prove that $z\in \mathcal{X}$. For each $z_{n}$ there exists $f_{n}\in S_{\Gamma}$ such that
$
z_{n}(t)=x_0+{\textrm {(P)}}\int_0^t k(t,s)f_{n}(s)ds.
$\\
By Theorem \ref{separable}, we can apply
\cite[Theorem V-13]{castaing valadier} (see also the remark at p. 147) in order to get
the compactness of $S_{\Gamma}$ with respect to the topology induced
by the tensor product $L^{\infty}([0,1])\otimes X^{\ast}$. This yields the compactness (and metrizability) of $S_{\Gamma}$ for the topology of convergence on the space of measurable simple functions taking a finite number of values in $X^*$ (described in \cite[ p. 176]{castaing valadier}). So, we
can extract a subsequence (not re-labelled) such that for each $t\in [0,1]$, ${\textrm {(P)}}\int_0^t
k(t,s)f_{n}(s) ds$ weakly converges to ${\textrm {(P)}}\int_0^t k(t,s)f(s) ds$. Now from the
uniform convergence of $z_{n}$ towards $z$ it follows that for all
$t\in [0,1]$,
$
z(t)=x_0+{\textrm {(P)}}\int_0^t k(t,s)f(s)ds
$ and so, $z\in \mathcal{X}$.\\
Consider now the multivalued operator $\Phi:\mathcal{X}\to c(\mathcal{X})$ defined by
$$
\Phi(x)=\left\{ y\in \mathcal{X} :  \  y(t)=x_0+{\textrm {(P)}}\int_0^t k(t,s)f(s)ds, f\in
S_{F(\cdot,x(\cdot))} \right\}.
$$
By Corollary \ref{cc} and  \cite[Theorem 2.9]{himmel}, its values are non-empty. As in \cite[Theorem VI-7]{castaing valadier} it follows that the values of $\Phi$ are closed (therefore, compact) and that $\Phi$ is upper semi-continuous.
By Kakutani-Ky Fan's Theorem, the operator $\Phi$ has a fixed point, which is a solution of the integral inclusion (1).
\hfill$\Box$\\

 Taking into account the notion of pseudo-derivative, we get also:

 \begin{cor}
Let $\Gamma:[0,1]\to ck(X)$ be Pettis-integrable and Bochner measurable and $F:[0,1]\times X\to ck(X)$ satisfy hypotheses 1)-3) of Theorem \ref{Volterra}.
Then the differential problem
\be\label{edp}
\qquad \qquad x'(t)\in F(t,x(t)),\; x(0)=x_0
\ee
 has at least one pseudo-differentiable solution $x\in C([0,1],X)$.
\end{cor}

The existence of solutions can be also gathered if we impose the integrability of the multifunction $\Gamma$. In this case, the solutions are continuous and a.e. differentiable in the usual meaning since, by   \cite[Theorem 3.2]{orlov ston}, in Fr\'{e}chet spaces, the class of indefinite Bochner-integrals coincides with the class of absolutely continuous and a.e. differentiable mappings.\\
\begin{cor}
Let $\Gamma:[0,1]\to ck(X)$ be integrable and $F:[0,1]\times X\to ck(X)$ satisfy hypotheses 1)-3) of Theorem \ref{Volterra}.
Then the differential problem (\ref{edp}) has at least one differentiable solution $x\in C([0,1],X)$.
\end{cor}

\normalsize



\end{document}